\documentclass{amsart}
\usepackage{amssymb}
\usepackage{amsmath}
\usepackage{amscd}
\usepackage{graphics}
\usepackage{latexsym}
\usepackage{amsrefs}
\usepackage{hyperref}

\def\R{\mathbb{R}}
\def\Z{\mathbb{Z}}

\newtheorem{theorem}{Theorem}[section]
\newtheorem{lemma}[theorem]{Lemma}

\newtheorem{corollary}[theorem]{Corollary}
\newtheorem{definition}[theorem]{Definition}
\newtheorem{remark}[theorem]{Remark}

\def\CP{\hbox{${\mathbb C}P^2$}}
\def\CPB{\hbox{$\overline{{\mathbb C}P^2}$}}

\def\w{\hbox{$\omega$}}
\def\s2s3{\hbox{$S^2\times S^3$}}

\title{On fillings of homotopy equivalent contact structures}
\author{Ahmet Beyaz}
\address{Department of Mathematics, Middle East Technical University, Ankara 06800 Turkey}
\email{beyaz@metu.edu.tr}
\thanks{}

\begin{document}
\begin{abstract}
This paper provides a topological method for filling contact structures on the connected sums of $S^2\times S^3$. Examples of nonsymplectomorphic strong fillings of homotopy equivalent contact structures with vanishing first Chern class on $\#_k S^2\times S^3$ $(k\geq2)$ are produced. 
\end{abstract}

\maketitle
\setcounter{section}{-1}

\section{Introduction} \label{introduction}

The study of four dimensional topology has seen great advances for the last $35$ years after Freedman (\cite{Freedman1982}) and Donaldson (\cite{Donaldson1983}). In the last several years, exotic smooth structures on rather small $4$-manifolds have been discovered (\cite{Fintushel2012}). This note exploits the symplectic structures on these simply connected $4$-manifolds with small second homology groups. The fillings of contact $5$-manifolds are distinguished by the lifts of symplectic surfaces inside the $4$-manifolds into the fillings. In general, fillings of contact manifolds may be used to extract information about the contact structure. In this paper the contact structures on the filled $5$-manifold are homotopy equivalent. However it is not clear whether these contact structures are isotopic or contactomorphic.

The paper consists of two sections. Section~\ref{preliminaries} is on preliminaries about contact structures and fillings of them. The second section includes the main theorems (Theorem~\ref{k2} and Theorem~\ref{k9}) and their proofs.

\section{Preliminaries} \label{preliminaries} This section reviews some definitons and facts about contact structures on $5$-manifolds. More information can be found in (\cite{Geiges2008}).

{\definition Let $N$ be a manifold of odd dimension $2n+1$. A contact structure is a maximally nonintegrable hyperplane field $\xi=kernel(\alpha) \subset TM$. The defining differential $1$-form $\alpha$ is required to satisfy $\alpha\wedge(d\alpha)^n>0$. Such a $1$-form $\alpha$ is called a contact form. The pair $(N,\xi)$ is called a contact manifold. A complex bundle structure $J$ on $\xi$ is called $\xi$-compatible if $J_p:\xi_p\rightarrow \xi_p$ is a $d\alpha$-compatible complex structure on $\xi_p$ for each point $p\in N$, where $\alpha$ is any contact form such that $\xi=kernel(\alpha)$. A $\xi$-compatible almost contact structure on a contact manifold $(N,\xi)$ is a complex structure on $\xi$ which is $\xi$-compatible.}\\

The condition $\alpha\wedge(d\alpha)^n>0$ means that the orientation of $N$ and the orientation imposed by the contact structure are same. $N$ is oriented by $\alpha\wedge(d\alpha)^n>0$ and $\xi$ is oriented by $(d\alpha)^{n-1}$. Note that $d\alpha$ gives a symplectic vector bundle structure to $\xi$.


{\definition Two contact manifolds $(N_0,\xi_0)$ and $(N_1,\xi_1)$ are said to be contactomorphic if there is an orientation preserving diffeomorphism $f:N_0\to N_1$ with $df(\xi_0)=\xi_1$, where $df:TN_0\to TN_1$ denotes the differential of $f$. If $\xi_i = kernel(\alpha_i)$, $i=0,1$, this is equivalent to saying that $\alpha_0$ and $f^*\alpha_1$ determine the same hyperplane field, and hence equivalent to the existence of a positive function $g: N_0\to\R^+$ such that $f^*\alpha_1 = g\alpha_0$.

Two contact structures $\xi_0$ and $\xi_1$ on a smooth manifold $N$ are said to be homotopy equivalent if their respective almost contact structures are homotopy equivalent. $\xi_0$ and $\xi_1$ are said to be isotopic if there is a smooth isotopy $\psi_t$ $(t\in [0,1])$ of $N$ such that $T\psi_t(\xi_0)=\xi_t$ for each $t\in[0,1]$. Equivalently, $\psi_t^*\alpha_t=\lambda_t\alpha_0$, where $\lambda_t:N\rightarrow \R^+$ is a suitable smooth family of smooth functions. This is equivalent to existence of a contactomorphism $f:(N,\xi_0) \to (N,\xi_1)$ which is isotopic to the identity.}\\

If two contact structures $\xi_0$ and $\xi_1$ on $N$ are isotopic then $(N,\xi_0)$ and $(N,\xi_1)$ are contactomorphic. Homotopy equivalence is much weaker than the isotopy and there may be many nonisotopic contact structures in a homotopy type.

A simply connected $5$-manifold $N$ admits an almost contact structure if and only if its integral Stiefel-Whitney class $W_3$ vanishes. Homotopy classes of almost contact structures are in one to one correspondence with integral lifts of $w_2(TN)$. The correspondence is given by associating to an almost contact structure its first Chern class (\cite{Geiges2008} p368).

{\definition 
A compact symplectic manifold $(M,\w)$ is called a strong (symplectic) filling of $(N,\xi)$ if $\partial M=N$ and there is a Liouville vector field $Y$ defined near  $\partial M$, pointing outwards along  $\partial M$, and satisfying $\xi= kernel(\w(Y,\cdot)|_{T M})$ (as cooriented contact structure). In this case we say that $(N,\xi)$ is the convex (or more precisely: $\w$-convex) boundary of $(M,\w)$.}\\

For contact manifolds of dimensions greater than three, $M$ is a strong filling of $N$ if and only if $\partial M=N$ as oriented manifolds and $\w|_{\xi}$ is in the conformal class of $d\alpha|_{\xi}$ (Theorem 5.1.5 of \cite{Geiges2008}). The boundary of a strong filling is said to be of contact type.

Let $(X,\omega)$ be a symplectic $4$-manifold and $e$ be a second cohomology class of $X$. Let's denote the $2$-disk bundle over $X$ with Euler class $e$ by $M_e$. Let $\pi$ be the projection map of the fibration. Any symplectic form on $M_e$ is locally $\pi^*\omega\oplus\omega_{a}$ where $\omega_{a}$ is the symplectic structure on the fiber for $a\in X$. Any two such symplectic forms agree on the zero section, therefore $M_e$ has a symplectic structure which is unique up to symplectomorphism by the symplectic neighborhood theorem (\cite{McDuff1998}). The next lemma is deduced from \cite{Etnyre1998}.

{\lemma \label{convex} Let $[\omega]$ be the cohomology class of the symplectic form. If $[\omega]e<0$ then the contact structure on the boundary is of contact type.}

\section{Nonsymplectomorphic Fillings of a Homotopy Type} \label{homotopy}

Assume that $X$ is a closed, simply connected, smooth $4$-manifold and $e \in H^{2}(X;\Z)$ is a primitive, characteristic class. If $X_e$ is the total space of the $S^1$-bundle over $X$ with Euler class $e$, then $X_e$ is diffeomorphic to $\#_{b_2(X)-1}S^2\times S^3$ (\cite{Duan2005}).



The pullback of the almost contact structure which is compatible with $d\alpha$ on the $5$-manifold is the pullback of the symplectic form $\w$ on the base $4$-manifold. The next lemma is relating the first Chern class of  the contact structure on the boundary and the symplectic structure on the base $4$-manifold.

{\lemma \label{contactSymplectic} The first Chern class of a compatible almost contact structure is the pullback of the first Chern class of the symplectic structure $\omega$ on $X$.}

\subsection{Fillings of $\#_k S^2\times S^3$ ($k\geq2$)} 

{\theorem \label{k2} In the homotopy equivalence class of contact structures on $\#_2 S^2\times S^3$ with the first Chern class equal to zero, there are contact structures which have nonsymplectomorphic strong fillings.}
\begin{proof} According to Fintushel and Stern (\cite{Fintushel2011,Fintushel2012}), there are infinitely many mutually nondiffeomorphic smooth manifolds which are homeomorphic to $\CP\#_2\CPB$, two of which carry symplectic structures. Let $(X_0,\omega_0)$ be $\CP\#_2\CPB$ and let $(X_1,\omega_1)$ be the symplectic $4$-manifold which is homeomorphic to $\CP\#_2\CPB$, but not diffeomorphic to it as given in \cite{Fintushel2012}. $X_0$ is not minimal and $X_1$ is minimal. 
Both of the manifolds have just two Seiberg-Witten basic classes which are plus and minus the canonical class. The first Chern class of $X_0$ is $3H-E_1-E_2$. On the other hand the canonical class of $X_1$ evaluates positive with the symplectic form, because it is a surface of general type (\cite{Fintushel2012} page 66). The first Chern class of $X_1$ is either $3H-E_1-E_2$ or $-(3H-E_1-E_2)$. 

Let $e_0$ be $-c_1(X_0)=-(3H-E_1-E_2)\in H^2(X_0;\Z)$ and $e_1$ be $c_1(X_1)$ in $H^2(X_1;\Z)$. For $j=0,1$, the boundary of $M_{e_j}$ is the circle bundle over $(X_j,\w_j)$ with Euler class $e_j$, and the smooth structures on the boundaries of  $M_{e_1}$ and $M_{e_2}$ are diffeomorphic to $\#_2 S^2\times S^3$ (\cite{Duan2005}). Since $X_j$ 
is simply connected, a part of the Gysin sequence for this circle bundle over $X_j$ is as shown below. 

\begin{equation} \label{Gysin}
0\rightarrow H^0(X_j;\Z)\overset{\cup e_j }{\rightarrow }H^{2}(X_j;\Z) \overset{\pi^*}{\rightarrow } H^{2}(\#_{b_2(X_j)-1}S^2\times S^3;\Z)\rightarrow0=H^1(X_j;\Z)
\end{equation}

The image of the the map $\cup e_j $ is generated by $e_j$ that is by plus or minus $c_1(X_j)$. By Lemma \ref{contactSymplectic}, the pullbacks of $c_1(\xi_0)$ and $c_1(\xi_1)$ on $\#_2 S^2\times S^3$ are the first Chern classes of the respective symplectic structures on $X_0$ and $X_1$. These classes are in the kernel of $\pi^*$, therefore the first Chern classes of the respective contact structures are zero.

The symplectic form $\w_0$ on $X_0$ couples negatively with $e_0$ (\cite{McDuff1994a}) and $\w_1\cdot e_1$ is less than zero by the discussion above. By Lemma~\ref{convex} the boundaries are of contact type. The boundaries are diffeomorphic to $\#_2 S^2\times S^3$ and the first Chern classes of the corresponding contact structures are zero. 

It remains to show that the symplectic structures are different. This is done by a count of $J$-holomorphic curves. For $j=0, 1$ the inclusion of $X_j$ into $M_{e_j}$ induces an injection of $H_2(X_j;\Z)$ into $H_2(M_{e_j};\Z)$. Let $\overline{E}$ be the image of $E\in H_2(X_j;\Z)$ in $H_2(M_{e_j};\Z)$. Remember $E\cdot E$ is $-1$ in $X_j$. Let $J_j$ be a generic almost complex structure which is compatible with the symplectic structure on $M_{e_j}$. Since $E$ is the class of an exceptional sphere in $X_0$, $E$ has an almost complex sphere representative in $X_0$. The image of this sphere under the inclusion map is a $J_0$-holomorphic sphere representative of $\overline{E}$ in $M_{e_j}$. Assume $\overline{E}$ has an $J_1$-holomorphic sphere $M_{e_1}$ that represents $\overline{E}$. Then $E$ would have a sphere representative in $X_1$. But this is not the case because $X_1$ is minimal. Therefore $M_{e_0}$ and $M_{e_1}$ are not symplectomorphic.
\end{proof}

{\remark It is not clear for the author that whether the symplectic structures on the $4$-manifolds $X_0$ and $X_1$ are related in any way. As a result of the reverse engineering which is applied to a model manifold, $X_1$ is known to be symplectic (\cite{Fintushel2007}, \cite{Fintushel2011}). In \cite{Fintushel2011} and \cite{Fintushel2012} this symplectic manifold $X_1$ is obtained from $X_0$ with a surgery on a single nullhomologous torus. But the latter operation does not involve the symplectic structures. So there is an ambiguity in the choice of $e_1$ in the proof.}

{\remark In \cite{Stipsicz2005c}, Stipsicz and Szabo note that Seiberg–Witten invariants can tell apart only at most finitely many symplectic structures on the topological manifold $\CP\#_k\CPB$ with $k\leq8$. Therefore this infinity result can not be extended to lower $k$ in an obvious way by the methods of this paper.
}

By using contact surgery and symplectic handlebody results of Meckert and Weinstein (\cite{Weinstein1991, Meckert1982}) one can say:

{\corollary For $k\geq2$, in the homotopy equivalence class of contact structures on $\#_k S^2\times S^3$ with the first Chern class equal to zero, there are contact structures which have nonsymplectomorphic strong fillings.}\\

\subsection{Fillings of $\#_k S^2\times S^3$ ($k\geq9$)} 

Dolgachev surfaces are elliptic surfaces that are homeomorphic to the elliptic surface $E(1)$ but not diffeomorphic to it. These manifolds are denoted by $E(1)_{p,q}$. $E(1)_{p,q}$ can be constructed from $E(1)$, which is diffeomorphic to $\CP\#_9\CPB$, by $p$ and $q$ logarithmic transformations where $gcd(p,q)=1$ and $p>q>1$. Considering the infinitely many different symplectic structures on these manifolds, one can conclude as follows.

{\theorem \label{k9} In the homotopy equivalence class of contact structures on $\#_9 S^2\times S^3$ with the first Chern class equal to zero, there are contact structures which have infinitely many nonsymplectomorphic strong fillings.}
\begin{proof} Assume that $p,q\in\Z$ such that $gcd(p,q)=1$ and $p>q>1$. If $F$ is the class of a generic fiber of the elliptic fibration on $E(1)_{p,q}$, then there is a homology class $A_{p,q}=\frac{F}{pq}$ in $H_2(E(1)_{p,q};\Z)$. The first Chern class of $E(1)_{p,q}$ is $-(pq-p-q)PD(A_{p,q})$, a negative multiple of Poincare dual of $A_{p,q}$. Let $e_{p,q}$ be $-PD(A_{p,q})$, which is primitive, and $M_{p,q}$ be the total space of the disk bundle over $E(1)_{p,q}$ with Euler class $e_{p,q}$. The boundary of $M_{p,q}$ is the circle bundle $X_{p,q}$ over $E(1)_{p,q}$ with Euler class $e_{p,q}$ and it is diffeomorphic to $\#_9 S^2\times S^3$ for all $p,q$. 

By Lemma \ref{contactSymplectic}, the pullbacks of the first Chern classes of the contact structures on $X_{p,q}$ the first Chern classes of the respective symplectic structures on $M_{p,q}$. The Gysin sequence for this circle bundle over $E(1)_{p,q}$ gives the first Chern classes of these contact structures are zero.

$E(1)_{p,q}$ are simply connected (proper) elliptic surfaces. According to the Kodaira-Enriques classification of complex surfaces, the symplectic form evaluates negatively on $c_1(E(1)_{p,q})$ and on $e_{p,q}$. By Lemma~\ref{convex}, for each $\{p,q\}$, the symplectic structure on the disk bundle $M_{p,q}$ is strong filling of its contact boundary.

Let $\overline{A}_{p,q}\in H_2(M_{p,q};\Z)$ be the pushforward of the class $A_{p,q}\in H_2(E(1)_{p,q};\Z)$. As explained by Ruan and Tian in \cite{Ruan1997} page 505, for the choice of a complex structure $J_{p,q}$, among all multiples of $mA_{p,q}$ (for $0<m<pq$), only $pA_{p,q}$ and $qA_{p,q}$ have connected $J_{pq}$-holomorphic torus representatives. Moreover this choice of complex structure is generic. This means, for any two distinct couples $\{p,q\}$ and $\{p',q'\}$, either $p\overline{A}_{p,q}$ and $q\overline{A}_{p,q}$ or $p'\overline{A}_{p',q'}$ and $q'\overline{A}_{p',q'}$ have connected complex torus representatives in the total space of the disk bundle. Therefore symplectic structures on $M_{p,q}$ and $M_{p',q'}$ are not symplectomorphic.
\end{proof}

{\corollary For $k\geq9$, in the homotopy equivalence class of contact structures on $\#_k S^2\times S^3$ with the first Chern class equal to zero, there are contact structures which have infinitely many nonsymplectomorphic strong fillings.}

\bibliography{library}

\end{document}